\documentclass{icmart}

\usepackage{color}

\contact[lerdos@ist.ac.at]{IST Austria, Am Campus 1, A-3400 Klosterneuburg, Austria}

\newtheorem{theorem}{Theorem}[section]
\newtheorem{corollary}[theorem]{Corollary}

\theoremstyle{definition}
\newtheorem{definition}[theorem]{Definition}

\renewcommand{\cal}{\mathcal}
\newcommand{\cU}{\cal{U}}
\newcommand{\cA}{\cal{A}}
\newcommand{\cB}{\cal{B}}
\newcommand{\cW}{\cal{W}}

\newcommand{\wt}{\widetilde}
\renewcommand{\P}{\mathbb{P}}
\newcommand{\E}{\mathbb{E}}
\newcommand{\R}{\mathbb{R}}

\newcommand{\N}{\mathbb{N}}

\newcommand{\e}{\varepsilon}

\newcommand{\pt}{\partial}
\newcommand{\rd}{{\rm d}}
\newcommand{\bR}{{\mathbb R}}
\newcommand{\bC}{{\mathbb C}}

\newcommand{\bZ}{{\mathbb Z}}
\newcommand{\non}{\nonumber}

\renewcommand{\Re}{\mbox{Re}}
\renewcommand{\Im}{\mbox{Im}}

\newcommand{\bx}{{\bf{x}}}

\newcommand{\bu}{{\bf{u}}}

\newcommand{\bla}{\mbox{\boldmath $\lambda$}}

\newcommand{\al}{\alpha}

\newcommand{\la}{\lambda}

\newcommand{\om}{{\omega}}
\newcommand{\si}{\sigma}

\newcommand{\cH}{{\mathcal H}}

\newcommand{\ov}{\overline}

\DeclareMathOperator{\tr}{Tr}

\newcommand{\llbracket}{[\![}
\newcommand{\rrbracket}{]\!]}

\newcommand{\be}{\begin{equation}}
\newcommand{\ee}{\end{equation}}

\newcommand{\nc}{\normalcolor}

\title[Random matrices, log-gases and H\"older regularity]{Random matrices, log-gases
and  H\"older regularity}

\author[L\'aszl\'o Erd{\H o}s]
{L\'aszl\'o Erd{\H o}s\thanks{Partially supported by SFB-TR 12 Grant of the German Research Council 
 and by ERC Advanced Grant, RANMAT 338804.}}

\begin{document}
\date{May 15, 2014}

\begin{abstract}
The Wigner-Dyson-Gaudin-Mehta conjecture asserts that 
the local eigenvalue statistics of large real and complex Hermitian matrices
with independent, identically distributed entries
are universal in a sense that  they depend only on 
the symmetry class of the matrix and otherwise are
 independent of the details of the distribution.
We present the recent solution to this half-century old
conjecture. We explain how stochastic tools, such as the 
Dyson Brownian motion, and  PDE ideas,
such as De Giorgi-Nash-Moser regularity theory, were combined in the solution.
 We also show related results for
log-gases that represent a universal model for
strongly correlated systems. Finally, in the spirit 
of  Wigner's original vision,
we discuss the extensions of these universality results
to  more realistic physical systems such as random band matrices.
\end{abstract}

\begin{classification}
Primary 15B52; Secondary 82B44.
\end{classification}

\begin{keywords}
De Giorgi-Nash-Moser parabolic regularity, Wigner-Dyson-Gaudin-Mehta universality, Dyson Brownian motion
\end{keywords}

\maketitle

\section{Introduction}

\bigskip

Large complex systems with many degrees of freedom  often exhibit remarkably
simple universal patterns. The Gauss law describes the fluctuations
of large sums of independent or weakly dependent random variables irrespective of
their distribution. The Poisson point process is the universal model for
many independent events in space or time. Both laws are ubiquitous in Nature
thanks to their large domain of attraction but they cannot accurately model
strong correlations. Can one find a universality for  correlated systems?

Since correlations appear in many forms, this seems an impossible task. 
Nevertheless this is exactly what E. Wigner
has accomplished when he discovered a universal pattern in the spectrum of
 heavy nuclei. Spectral measurement data for various nuclei clearly show
that the density of energy levels depends on the actual nucleus.
But Wigner asked a different question: he looked 
at the energy {\it gaps}, i.e.
the {\it difference} of consecutive energy levels.
He discovered that their statistics, after rescaling with the local density, 
showed a very similar pattern for different nuclei. 

Wigner's revolutionary insight was that
this coincidence does not stem from some particular
property of the specific physical system but it
has a profound mathematical origin. General quantum mechanics
postulates that energy levels are eigenvalues of 
a certain hermitian matrix (or operator) $H= (h_{ij})$, the {\it Hamiltonian} of the system. 
The matrix elements $h_{ij}$  represent quantum transition
rates between two states labelled by $i$ and $j$. While $h_{ij}$'s are specific to the system,
the gap statistics largely depend only on the basic symmetry class of $H$,
as long as $h_{ij}$'s are chosen somewhat generically.

To illustrate this mechanism, consider a $2\times 2$ hermitian matrix
$$
   H = \begin{pmatrix} a & b\cr \bar b & d \end{pmatrix}, \qquad a,d \in \bR, \; b\in \bC.
$$
The difference (or gap) of the two eigenvalues is $\lambda_2-\lambda_1 =\big[ (a-d)^2 + 4|b|^2 \big]^{1/2}$.
If the matrix elements are drawn independently from some continuous distribution, then the probability
that the gap is very small; 
$$
   \P (|\lambda_2-\lambda_1|\leq \e), \qquad \e\ll 1,
$$
is of order $\e^2$ for real symmetric matrices ($b\in \bR$) and it is
of order $\e^3$ for complex hermitian matrices ($b\in \bC$). The exponent
of $\e$ is thus determined by the symmetry class of $H$.

Very surprisingly, for large $N\times N$ matrices
the {\it entire distribution} of the gap becomes universal as $N\to\infty$ and 
not only its asymptotics in the $\e\ll 1$ regime.
Based upon a
more precise calculation with Gaussian matrix elements, Wigner
 predicted that this  universal law is  given by a simple
formula (called the {\it Wigner surmise}). For the real symmetric case it is 
\be\label{surmise}
   \P\Big( \wt \lambda_j-\wt \lambda_{j-1} 
= s+\rd s\Big) \approx 
\frac{\pi s}{2}
\exp\big( -  \frac{\pi}{4} s^2\big)\rd s,
\ee
where $\wt \lambda_j = \varrho \lambda_j$ denote the 
eigenvalues $\lambda_j$ rescaled by the density of eigenvalues $\varrho$ near $\lambda_j$. 
This law is characteristically different from the gap
distribution of the Poisson point process which is the exponential distribution, $e^{-s}\rd s$.
The prefactor $s$ in \eqref{surmise} indicates a {\it level repulsion}
 for the point process $\wt \lambda_j$, in particular the eigenvalues are strongly correlated (eigenvalues are often called (energy) levels in 
random matrix theory).
Similar formulas hold for the joint statistics of several consecutive gaps.

Comparing measurement data from various experiments, Wigner concluded
that the energy gap distribution of complicated quantum systems is essentially universal;
it depends only  on the basic symmetries of model
(such as time-reversal invariance). This thesis 
has never been rigorously proved for any realistic 
physical system but experimental data and extensive 
numerics leave no doubt on its correctness (see \cite{M} for an overview).

Once universality is expected, explicit formulas for the 
statistics can be computed from the most convenient model
within the universality class. The simplest representatives
of these universality classes are $N\times N$ random matrices
with independent (up to symmetry), identically distributed Gaussian entries.
These are called the {\it Gaussian orthogonal ensemble (GOE)}
and the {\it Gaussian unitary ensemble (GUE)} in case of
real symmetric and complex Hermitian matrices, respectively.
Wigner's bold vision was to neglect all details
of the actual Hamiltonian operator  and replace it with a 
large random Gaussian matrix of the same symmetry class. As far as the
gap statistics are concerned, this simple-minded model
very accurately reproduced the behavior of large complex quantum systems!

 Since Wigner's discovery random matrix statistics are found 
everywhere in physics and beyond, wherever nontrivial correlations  prevail.
Random matrix theory (RMT) is present in chaotic quantum systems
in physics, in principal component analysis in statistics,
in communication theory and even in number theory. In particular,
the zeros of the Riemann zeta function on the critical line are
expected to follow RMT statistics due to a spectacular result of Montgomery \cite{Mont}.

In retrospect, Wigner's idea should have received even more attention.
For centuries, the primary territory of probability theory was to model uncorrelated
or weakly correlated systems.  The surprising ubiquity of random matrix statistics is
a strong evidence that it plays a similar fundamental role for correlated
systems as Gaussian distribution and Poisson point process play for
uncorrelated systems. RMT seems to provide essentially
 the only universal and generally   computable pattern  for complicated
 correlated systems.

A few years after Wigner's seminal paper \cite{W},
Gaudin \cite{Gau}  discovered another remarkable
property of this new point process:  the correlation
functions have an exact determinantal structure, at least
if the distributions of the matrix elements are Gaussian.
The algebraic identities within the determinantal form
 opened up the route  to obtain explicit formulas
for local correlation functions. For example, in the
complex Hermitian case (GUE)
the $n$-point correlation function $p^{(n)}$ of the
rescaled eigenvalues $\wt \la_i$ in the bulk  is given by the determinant of the celebrated sine-kernel:
\be\label{sine}
   p^{(n)}(\wt\la_1, \wt\la_2, \ldots, \wt\la_n) = \mbox{det}\big[ K(\wt\la_i, \wt\la_j)\big]_{i,j=1}^n, 
\qquad K(x, y) := \frac{\sin \pi(x-y)}{\pi (x-y)}.
\ee
(The same determinantal expression with a different but closely related kernel function $K$ holds
for the real symmetric case.)
As a consequence, the gap distribution  is given by a Fredholm determinant
involving Hermite polynomials. In fact,  Hermite polynomials were 
first  introduced in the context of random matrices by Mehta and Gaudin \cite{MG} earlier.
 Dyson and Mehta \cite{M2, Dy1, Dy2}
have later extended this exact calculation 
to correlation functions and to other symmetry classes.
When compared with the exact formula,
the Wigner surmise \eqref{surmise}, based upon a simple $2\times 2$
matrix model, turned out to be quite accurate.
While the determinantal structure is present only in
Gaussian Wigner matrices, the paradigm of spectral universality
predicts that the formulas for the local eigenvalue statistics
obtained in the Gaussian case hold for general distributions as well.

\section{Random matrix ensembles and log-gases}

We consider $N\times N$ hermitian matrices $H$ with matrix elements
having mean zero and variance $1/N$, i.e.
\be
   \E\, h_{ij} = 0,   \quad \E |h_{ij}|^2= \frac{1}{N}
 \qquad i,j =1,2,\ldots, N.
\label{centered}
\ee
The matrix elements  $h_{ij}$ are
real or  complex  independent random variables 
subject to the symmetry constraint $h_{ij}= \ov h_{ji}$.   These ensembles of random matrices 
are called {\it (standard) Wigner matrices.} The normalization \eqref{centered} is
introduced for definiteness.

An important special case of Wigner matrices is the Gaussian case (GOE or GUE), 
when $h_{ij}$'s have Gaussian distribution. In this case the matrix ensemble
can also be given by the probability law
\be\label{inv1}
   P (H) \rd H = Z^{-1} e^{-\frac{\beta}{4} N\tr H^2} \rd H,
\ee
where $\rd H = \prod_{i<j} \rd h_{ij} \rd\bar h_{ij}\prod_i \rd h_{ii}$ 
is the standard Lebesgue measure on  real symmetric or complex hermitian $N\times N$ matrices
and $Z=Z_N$ is the normalization. 
The parameter $\beta$ is chosen to be $\beta=1$ for GOE and $\beta=2$ for GUE
to ensure the normalization \eqref{centered}.

The representation \eqref{inv1}
shows that the Gaussian ensembles enjoy an invariance property; the distribution
$P(H)$ is invariant under a base transformation, $H\to UHU^*$, where $U$ is
orthogonal (in case of GOE) or unitary (in case of GUE). In fact, invariance property is
not restricted to the Gaussian case; one may directly generalize \eqref{inv1} to
\be\label{inv2}
   P (H) \rd H = Z^{-1} e^{-\frac{\beta}{2} N\tr V(H)} \rd H,
\ee
where $V:\bR\to\bR$ is an arbitrary function with sufficient
growth at infinity to ensure the normalizability of the measure. 
The ensembles of the form \eqref{inv2} are called {\it invariant matrix ensembles.}

Wigner ensembles and invariant ensembles represent two natural but quite different 
ways to equip the space of $N\times N$ matrices with a probability measure. 
These two families are essentially disjoint; only the Gaussian ensembles
belong to their intersection.

Let $\bla = (\lambda_1, \lambda_2, \ldots , \lambda_N)$ denote the eigenvalues of $H$
in increasing order.
Since eigenvalues are complicated functions of the matrix elements, there is no
explicit formula to express the probability distribution of $\bla$ induced by
a general Wigner ensemble. However, quite remarkably, for invariant ensembles \eqref{inv2}
the joint probability density of the eigenvalues is explicitly given by
\be\label{pinv}
  \mu_{\beta, V}^{(N)} (\bla)= 
C \prod_{1\le i<j\le N} (\lambda_j-\lambda_i)^\beta 
\prod_{j=1}^N e^{-\frac{\beta}{2} N V(\lambda_j)}
\ee
with a normalization constant $C$. This formula may directly be obtained from
\eqref{inv2} by diagonalizing $H=U\Lambda U^*$ and integrating out the
matrix of eigenvectors $U \in O(N)$ or $U\in U(N)$ with respect to the Haar measure.

From statistical physics point of view, we may consider the 
distribution \eqref{pinv} as a Gibbs measure for a gas of 
$N$ point particles on $\bR$. We may write
\be\label{loggas}
  \mu_{\beta, V}^{(N)} (\bla) = C\, e^{-\beta N\cH(\bla)}, \qquad  \cH(\bla):=  \sum_{k=1}^N  \frac{1}{2}V(\lambda_k)- 
\frac{1}{N} \sum_{1\leq i<j\leq N}\log (\lambda_j-\lambda_i),
\ee
where, according to  the Gibbs formalism, $\cH(\bla)$ is the Hamiltonian (energy function) of the gas
and the parameter $\beta$ plays the role of the inverse temperature. 
The Vandermonde determinant in \eqref{pinv}  translates into a logarithmic pair
interaction between the particles.
We may  completely ignore the original random matrix ensemble behind \eqref{pinv}
and consider \eqref{loggas} more generally for any parameter $\beta>0$, not only for the specific
values $\beta=1, 2$. The Gibbs measure \eqref{loggas} is often called {\it $\beta$-log-gas}
or {\it $\beta$-ensemble}.

Eigenvalue distributions of Wigner ensembles and
  $\beta$-log-gases are quite different mathematical
entities despite their connection via \eqref{pinv} in the special Gaussian case,
$V(\lambda)=\frac{1}{2}\lambda^2$ and $\beta=1,2$. 
Wigner ensembles are parametrized by the value $\beta=1, 2$ and by
the distribution of the single matrix elements,
while log-gases are parametrized by $\beta$ and the potential function $V$.
The central thesis of universality asserts that the gap statistics of both families of ensembles
depend only on the parameter $\beta$ and are otherwise
 independent of any other details of the models.

For Wigner matrices this thesis is generally referred to as the {\it universality conjecture 
 of  random matrices}  and we will call it the {\it Wigner-Dyson-Gaudin-Mehta conjecture}.
 It was first formulated in Mehta's treatise on
random matrices \cite{M} in 1967 and has remained a  key question
in the subject ever since. In this article we
 review the recent progress that has led 
to the proof of this conjecture and the analogous conjecture
for log gases. For more
details, the reader is referred to the lecture notes \cite{ECDM}.

\section{Random band matrices and Anderson model}

As mentioned in the introduction,
Wigner's vision extends the thesis of universality far beyond the models we just introduced. We now present
an extension that was an important source of motivation in the development of the subject.

Viewed as a quantum mechanical Hamilton operator, a Wigner matrix $H$ represents a {\it mean-field system}; the 
quantum transition rates $h_{ij}$ between any two quantum states, labelled by $i$ and $j$, are
comparable in size. The quantum states of more realistic physical models have a spatial structure
and typically quantum transition occurs between nearby states only. 

The spatial structure is essential to understand the  {\it metal-insulator transition} which 
is the fundamental phase transition of disordered quantum systems modelled by a 
random Hamilton operator  $H$. According to the physical theory, 
in the {\it metallic phase} the eigenfunctions are delocalized, the quantum time evolution $e^{itH}$
is diffusive and the local eigenvalue statistics coincide with the ones
 from the GUE/GOE random matrix theory
\eqref{sine}. The {\it localization length}, which is
the  characteristic lengthscale of the physically relevant quantities
(such as eigenfunctions or propagators), is practically infinite. In contrast, in the  {\it insulator phase},
the eigenfunctions are localized with a localization length $\ell$ independent of the system size,
the time evolution remains bounded for all times and the local eigenvalue statistics are Poisson.
In the mathematics literature these two phases are usually called {\it delocalized} and {\it localized} regimes, respectively,
and they are primarily characterized by the spectral type (absolutely continuous vs. pure point)
of the corresponding infinite volume operator.

The basic model for the metal-insulator transition is the celebrated
 Anderson model in solid state physics \cite{A}.
The Anderson Hamiltonian is given by $-\Delta + V(x)$ on the Hilbert space $\ell^2(\bZ^d)$, where
$\Delta$ is the lattice Laplacian and $V(x)$ is a real valued random potential field such
that $\{ V(x)\;: \; x\in \bZ^d\}$ are independent and 
identically distributed centered random variables with variance $\sigma^2:=\E |V(x)|^2$.
The Anderson model has been extensively studied mathematically. In nutshell, the
high disorder regime is relatively well understood since the seminal
work of Fr\"ohlich and Spencer \cite{FS} for localization (an alternative proof is given by
Aizenman and Molchanov \cite{AM}), complemented by the work of Minami \cite{Min} 
proving the local Poissonian spectral statistics.  In contrast, in the low disorder regime, starting from three 
spatial dimension and away from the spectral edges, the eigenfunctions
are conjectured to be delocalized but no rigorous proof exists  ({\it extended states conjecture}).

Random band matrices  are another popular model for the metal-insulator transition \cite{Spe}.
For definiteness, let the state space be
 a finite box $\Lambda:=[1, L]^d \subset \bZ^d$ of the $d$-dimensional integer lattice equipped with periodic boundary condition.
We consider hermitian matrices $H = (h_{ij})_{i,j\in \Lambda}$ whose rows and columns are labelled
 by the elements of $\Lambda$ and whose matrix
elements are independent.  Given a parameter $W\le L/2$, called the {\it  band width}, we assume
that the matrix elements $h_{ij}$ vanish beyond a distance $|i-j|\ge W$, i.e. we replace  
\eqref{centered} with the condition
\be
   \E\, h_{ij} = 0,   \;\; \forall i,j\in \Lambda; \quad\mbox{and}\quad  h_{ij}= 0 \quad \mbox{for}\quad |i-j|\ge W. 
\label{centered2}
\ee
($|\cdot |$ denotes the periodic distance on $\Lambda$). These are called {\it random band matrices}.
We often assume a translation invariant profile for the variances, i.e. that
\be\label{profile}
 \sigma_{ij}^2: = \E |h_{ij}|^2 = \frac{1}{W^d} f\Big( \frac{|i-j|}{W}\Big)
\ee
with some compactly supported function $f\ge 0$ on $\bR^d$ with $\int_{\bR^d} f =1$.
Notice that the normalization is chosen such that
\be\label{rowsum}
\sum_{j\in \Lambda}\sigma_{ij}^2=1, \qquad \forall i\in\Lambda.
\ee

If the band width is maximal,  $W= L/2$,,
 and $f$ is constant on $[-\frac{1}{2}, \frac{1}{2}]^d$, then we recover the Wigner matrices \eqref{centered}.
Wigner matrices are always in the delocalized regime as it was shown that all eigenfunctions are
extended with very high probability \cite{ESY2}.
The other extreme is when $W$ remains bounded even as the matrix size $|\Lambda|=L^d$ goes to infinity.
This system behaves very similarly to the  Anderson model.
In particular, in $d=1$  it exhibits Anderson localization  even if $W$ grows slowly with $L$ as $W\ll L^{1/8}$ \cite{Sche}.
 Therefore random band matrices with an intermediate band width, $1\ll W \ll L$, serve as a 
model to study the metal-insulator transition. The fundamental conjecture in $d=1$ 
 is that the transition occurs at $W=L^{1/2}$. 
This conjecture is supported by supersymmetric (SUSY) functional integration techniques \cite{Fy}
which are intriguingly elegant but notoriously hard to justify with full mathematical rigour.
Nevertheless, very recently  sine-kernel local statistics \eqref{sine} were proven 
for a  Gaussian band matrix with a specifically chosen block structure \cite{Scherb}
using SUSY approach. The details have been worked out for $W\ge L^{1-\e}$ with some small $\e>0$.
In a related problem (correlation function of the characteristic polynomial of $H$ at
two different energies) the result even holds down to the critical band width $W\ge L^{1/2+\e}$ \cite{Scherb2},
but still only for a specific block structure and Gaussian distribution. 

For more general band matrices the universality of the
local statistics have not yet been proven, but
it was shown in $d=1$ that the localization length is at least $W^{5/4}$, indicating 
band matrices with band width at least $W\gg L^{4/5}$ are in the delocalized regime
\cite{EKYY3}.

\section{Universality on three levels}

We consider an ensemble of  $N$ (unordered) random points $\bla=(\lambda_1, \lambda_2, \ldots, \lambda_N)$ 
on the real line, either given by eigenvalues of hermitian random
matrices or points of a log-gas. We always choose the normalization
such that all points lie in a bounded interval, independent of $N$, with a very high
probability. The typical spacing between the points is therefore of order $1/N$.

The statistics of $\bla$ are characterized by the $n$-point functions $p^{(n)}_N$. They are  defined by the
following relation that holds for any function $O$ of $n$ variables:
\begin{align}\label{def:corr}
  \E {N\choose n}^{-1} \sum^*
O(\lambda_{i_1},  & \lambda_{i_2}, \ldots, \lambda_{i_n}) 
\\ \nonumber
&  = \int_{\bR^n} p^{(n)}_N(x_1, x_2, \ldots , x_n) O(x_1, \ldots , x_n)\rd x_1 \ldots \rd x_n.
\end{align}
Here the star indicates that the summation runs over all $n$-tuples of distinct integers, $(i_1, i_2, \ldots, i_n)$ with
$1\le i_j\le N$. The correlation function for $n=1$ is called the {\it density}. Typically
we fix $n$ and consider the limit of the correlation functions $p^{(n)}_N$ as $N\to\infty$
to obtain the limiting statistics.

We may consider the limiting statistics of the points on three scales. For definiteness we 
 illustrate these scales for Wigner matrices; similar results hold for the log-gases and 
for random band matrices, but the latter only under more restrictive conditions.

\subsection{Macroscopic scale}

The largest  scale
corresponds to observable functions $O$ in \eqref{def:corr} that are unscaled with $N$. 
For Wigner matrices \eqref{centered} the limiting density is given by the celebrated semicircle law  \cite{W}
\be\label{def:sc}
   \varrho_{sc}(x): =\frac{1}{2\pi}\sqrt{ (4-x^2)_+}
\ee
in the form
of a weak limit:
\be\label{sclaw}
   \E \frac{1}{N}\sum_i O(\lambda_i) = \int_\bR p^{(1)}_N (x) O(x)\rd x \to \int_\bR \varrho_{sc}(x) O(x)\rd x, \quad
\mbox{as $N\to\infty$,}
\ee
that holds for any continuous, compactly supported function $O$. In fact, the semicircle law also holds not only in
expectation but also as a  convergence in probability for the empirical density:
\be\label{scprob}
  \P \Big(  \Big|  \frac{1}{N}\sum_i O(\lambda_i) -\int_\bR \varrho_{sc}(x) O(x)\rd x\Big|\ge\e\Big) \to 0
\ee
for any $\e>0$ as $N\to\infty$. 

These results are the simplest form of
spectral universality; they assert that the eigenvalue density on macroscopic scales is independent
of the specific distribution of the matrix elements. In fact,  this result also holds for 
{\it generalized Wigner matrices} whose matrix elements are still centered and independent,
but their distributions may vary. The semicircle law \eqref{sclaw} holds as long as the row sums of the variances
is constant, i.e.  
\be\label{rowsum1}
   \sum_j \sigma_{ij}^2 =1, \qquad \sigma_{ij}: = \E |h_{ij}|^2,
\ee
for any $i$. If \eqref{rowsum1} does not hold but the variances have a macroscopic profile
in a sense that
$$
   \sigma_{ij}^2 = S\big( \frac{i}{N}, \frac{j}{N}\big)
$$ 
with some fixed function $S$ on $[0,1]\times [0,1]$, then the limiting density
still exists and can be computed from $S$,
 but it is not given by the semicircle law any more \cite{AEK}. These
results show that the limiting density is determined by variances of the matrix
elements alone and not by their full distribution.

\subsection{Mesoscopic scales}

We now consider an $N$-dependent scaling parameter $\eta=\eta_N>0$ and a fixed point $E$ in the
support of the limiting density, $|E|<2$
(real numbers $E$ in the context of location in the spectrum  are often called {\it energy} 
due to the physical meaning of the spectrum).
The regime $1/N\ll\eta\ll 1$ corresponds to mesoscopic scales;
on these scales the fluctuation of the empirical density around the semicircle density profile is 
still negligible, but
the effects of individual points are not yet visible.

We rescale the
observable around $E$ in a window of size $\eta$ and consider
\be\label{sclawlocal}
   \E \frac{1}{N\eta}\sum_i O\Big( \frac{\lambda_i-E}{\eta}\Big)
 = \int_\bR p^{(1)}_N (E+x\eta) O(x)\rd x .
\ee
If $\eta\to 0$ as $N\to\infty$, then formally \eqref{sclaw} would indicate that the limit of \eqref{sclawlocal}
is $\varrho_{sc}(E) \int O(x) \rd x$. This is indeed correct,  with some technical assumptions
even in the stronger sense \eqref{scprob},
as long as  $1/N\ll\eta \ll 1$. This is called the {\it local semicircle law}
in the bulk of the spectrum.
The first result  down to the optimal scale $\eta\gg 1/N$ (modulo $\log N$ factors) was given in
 \cite{ESY2} followed by several improvements and generalizations, see \cite{EKY3} for a summary. 
In particular, local semicircle law has also been extended to the spectral edge, $|E|=2$, where the optimal
scale is $\eta\gg N^{-2/3}$ reflecting the fact that the eigenvalue spacing near the edge is of
order $N^{-2/3}$.

Local semicircle laws imply, among others, that
the points $\lambda_j$ are very close to their {\it classical location}
denoted by $\gamma_j$ and
defined as the $j$-th quantile of the limiting density:
\be\label{def:gamma}
     \int_{-\infty}^{\gamma_j}\varrho_{sc}(x)\rd x = \frac{j}{N}.
\ee
More precisely, we have for any $j$ (including the extreme eigenvalues near the 
spectral edge) that
\be
   |\lambda_j -\gamma_j| \lesssim |\gamma_{j+1}-\gamma_j| 
\label{rig}
\ee
with a very high probability, where $\lesssim$ indicates logarithmic factors
\cite{EYYrigi}.
The property \eqref{rig} is called {\it rigidity} and it asserts that
the fluctuation of the points is essentially on the scale of the local
gap $|\gamma_{j+1}-\gamma_j|$. In particular, for points in the bulk spectrum,
their fluctuation is only slightly larger than $1/N$. 

Local semicircle law also holds 
for random band matrices with \eqref{profile}, however the local density is
controlled only down to scales $\eta\gg W^{-1}$,
 see \cite{EKY3} for a summary and also \cite{Sod}.
The regime $\eta \ll W^{-1}$ is mathematically unexplored
and there is no optimal rigidity result.

 While the density on mesoscopic scales behaves exactly as on the macroscopic scale, 
the density-density correlation exhibits a new
universality. For  two random variables, $X,Y$,
let $\langle X; Y\rangle = \E XY - \E X \, \E Y$ denote their covariance.
Consider two energies $E_2\ge E_1$ and a scale $\eta$ such that $N^{-1/7}\ll\eta \ll E_2-E_1\ll 1$. Then  for Wigner matrices
the covariance decays with a universal power-law \cite{BKh1, EK3, EK4}
\be\label{AS}
   \Bigg\langle \frac{1}{N\eta}\sum_i O\Big( \frac{\lambda_i-E_1}{\eta}\Big)\, ; \,  
\frac{1}{N\eta}\sum_i O\Big( \frac{\lambda_i-E_2}{\eta}\Big) \Bigg\rangle \sim - \big[ N(E_2-E_1)\big]^{-2}
\ee
(for Gaussian case the result extends to $\eta\gg 1/N$ \cite{BKh}).
Higher order moments  satisfy the
Wick theorem asymptotically, i.e. the local densities at different energies converge to a
  Gaussian variables with a non-trivial covariance structure \cite{EK3, EK4}.

Similar result holds for band matrices with
\eqref{profile} in $d$ dimensions, but the power law decay in \eqref{AS} undergoes a phase transition. 
For $W^{-d/7} \ll \eta\ll (W/L)^2$ the asymptotics \eqref{AS} holds with  
the mean-field exponent $-2$, while for $(W/L)^2\ll\eta\ll 1$ the power
in the right hand side becomes $-2+\frac{d}{2}$ for $d=1,2,3$
and it is logarithmic for $d=4$. In higher dimensions, $d\ge5$, 
the universality breaks down.
 This feature is closely related to the quantum diffusion phenomenon for
the unitary time evolution \cite{EK, EK2}. In the physics literature these asymptotics are called the Altshuler-Shklovskii
formulas and recently they have been rigorosly proved \cite{EK3, EK4}.

\subsection{Microscopic scale}

The most intriguing regime for universality is the microscopic scale 
where the scaling parameter $\eta$ in the observable is chosen comparable
with the typical local eigenvalue spacing. In particular,
individual eigenvalues are observed. This is the regime
for the gap distribution in Wigner's surmise, and
the original conjecture of Mehta \cite{M} on random matrix
universality  also pertains microscopic scales.

Before we formulate the precise results, we make two remarks
to explain why there will be different universality theorems.

First, for the local statistics
we need to distinguish  the bulk spectrum where $\eta\sim 1/N$
and the edge spectrum where $\eta\sim N^{-2/3}$.
Not only the scaling  but also the explicit
formulas are different in these two regimes.  The correlation functions
are asymptotically  determinantal (Pfaffian) in both cases, but
in the bulk they are given by the Dyson sine kernel \eqref{sine}
and its real symmetric counterpart, while at the edge
they are given by the Airy kernel \cite{TW, TW2}.
In all cases the explicit formulas have been computed 
in the corresponding Gaussian model which is computationally
the most accessible case via orthogonal polynomials.
The significance of orthogonal polynomials in random matrices has first been
realized by Gaudin, Mehta and Dyson \cite{Gau, MG, Dy2}. Their approach was
later generalized and combined
with the Riemann-Hilbert method to yield explicit asymptotic calculations
for broader classes of invariant ensembles, see \cite{BI, De1,  DKMVZ1, DKMVZ2,FIK, Lub, PS:97, PS} 
for the extensive literature in the $\beta=2$ case and  \cite{DG1, KS, Sch} for  
the more complicated $\beta =1, 4$ case. Our universality results show
that the local statistics for a general Wigner matrix or invariant ensemble
(or even more generally a $\beta$-log-gas) coincide with those
of the corresponding Gaussian model. Therefore all explicit asymptotic
calculations apart from the simplest Gaussian case become redundant.

Second, there is a subtle difference between the universality of 
$n$-point local correlation functions around a {\it fixed energy} $E$ and the universality
of $n$ consecutive points $\lambda_{j+1}, \lambda_{j+2}, \ldots \lambda_{j+n}$ for some
{\it fixed label} $j$.
The former asks for identifying the limit
\begin{align}\label{univlocal}
   \E \frac{1}{(N\eta)^n}\sum^*  & O\Big( \frac{\lambda_{i_1}-E}{\eta}, \; \frac{\lambda_{i_2}-E}{\eta}, \;
\ldots,  \frac{\lambda_{i_n}-E}{\eta} \Big) \\
\nonumber
& = \int_{\bR^n} p^{(n)}_N (E+x_1\eta, E+x_2\eta, \ldots , E+x_n\eta) O(x_1, \ldots , x_n)\rd x_1\ldots \rd x_n 
\end{align}
for any smooth, compactly supported observable $O$, i.e. identifying the weak 
limit of the rescaled correlation functions $p^{(n)}_N (E+x_1\eta, E+x_2\eta, \ldots )$ in the variables $x_1, \ldots , x_n$,
The latter asks for the joint distribution of $\lambda_{j+1}, \lambda_{j+2}, \ldots \lambda_{j+n}$ with an appropriate rescaling.

The rigidity \eqref{rig} locates the $j$-th eigenvalue $\lambda_j$ 
around a fixed energy $E=\gamma_j$ but only with a precision slightly larger than $1/N$. In 
fact, for the Gaussian ensembles  it is known \cite{Gus, ORour} that $\lambda_j-\gamma_j$ is Gaussian and it fluctuates on scale
$\sqrt{\log N}/N$ therefore there is no direct translation between the two types of universality results.
In particular, the universality of $n$ consecutive gaps which was
originally advocated by Wigner,  i.e. the limit of
\begin{align}\label{gapuniv}
   \E \,  O\Big(  \frac{\lambda_{j+1}-\lambda_{j}}{\eta}, & \;  \frac{\lambda_{j+2}-\lambda_{j+1}}{\eta},  \; 
\ldots,  \frac{\lambda_{j+n}-\lambda_{j+n-1}}{\eta}
\Big) \\
\nonumber
& = \int_{\bR^n} g^{(j)}_N (x_1,  x_2, \ldots , x_n) O(x_1, \ldots , x_n)\rd x_1\ldots \rd x_n,
\end{align}
with the natural scaling $\eta = 1/N$,
cannot be concluded from the fixed energy universality \eqref{univlocal}.

Given the historical importance of the Wigner surmise, it is somewhat surprising
that gap universality with a fixed label did not receive much attention until very recently. 
The first results on  the Wigner-Dyson-Gaudin-Mehta
 universality  proved \eqref{univlocal} in the sense of {\it average energy}, i.e.
after taking average in the parameter $E$ in a small interval of size $N^{-1+\e}$. 
Since $N^{-1+\e}$ is above the rigidity scale, average energy universality easily
implies {\it average label  gap universality}, i.e. the averaged version
of \eqref{gapuniv} after averaging
the label $j$ in an interval of size $N^\e$.

Our more recent understanding shows that there is a profound difference between
the weaker {\it ``averaged''} results versus the stronger {\it ``fixed''} ones. 
Obviously, {\it ``fixed''} results are necessary for the precise statistics of individual points
hence for  fully characterizing the limiting process. 
At first sight, removing the local averaging may only seem a fine technical point;
it merely requires to exclude the pathological case that a certain energy $E$ (or a certain label $j$)
might behave very differently than a typical one. Physicists have never worried
about this situation since there is no apparent reason for such pathology
(in fact Mehta's original version of the conjecture
did not specify the precise formulation of universality). 
Mathematically, however, it turned out surprisingly involved to exclude
 the worst case scenarios and  we needed to develop a completely new approach.
Finally, we point out that, unlike their averaged counterparts, the fixed energy
and the fixed label results are not equivalent, in fact each required 
a separate proof.

\section{Universality of local statistics: the main results}

\subsection{Wigner ensembles}

Our main results hold for a larger class of ensembles than the standard Wigner matrices, which
we will call {\it generalized Wigner matrices}.

\begin{definition}\label{def:genwig}(\cite{EYY})
 The real symmetric or complex Hermitian
matrix ensemble $H$ with centred and independent matrix elements $h_{ij}=\ov{h}_{ji}$, $i\le j$,
is called {\it generalized Wigner matrix} if 
the variances  $\si^2_{ij}= \E |h_{ij}|^2$ satisfy:
\begin{description}
\item[(A)] For any $j$ fixed
\be
   \sum_{i=1}^N \si^2_{ij} = 1 \, .
\label{sum}
\ee

\item[(B)]   There exist two positive constants, $C_1$ and $C_2$,
independent of $N$ such that
\be\label{1.3C}
\frac{C_1}{N} \le \si^2_{ij}\leq \frac{C_2}{N}.
\ee
\end{description} 
For Hermitian ensembles, we  additionally  require
 that for each $i,j$ the $2\times 2$ covariance matrix is bounded by $C/N$
in  matrix sense, i.e.
$$
\Sigma_{ij} :\;=\; 
\begin{pmatrix} \E (\Re h_{ij})^2 & \E (\Re h_{ij})(\Im h_{ij}) \\
  \E (\Re h_{ij})(\Im h_{ij}) & \E ( \Im h_{ij})^2
\end{pmatrix} \geq\; \frac{C}{N}.
$$
\end{definition}
The following theorem settles the average energy version of the Wigner-Dyson-Gaudin-Mehta
conjecture for generalized Wigner matrices. It is formulated under the weakest moment
assumptions. The same result under somewhat more restrictive assumptions
were already obtained in \cite{ESY4, ESYY}; 
 see also \cite{TV} for the complex hermitian case
and for a quite restricted class of real symmetric matrices. 
More details on the history can be found in \cite{ECDM}.

\begin{theorem}[Universality with averaged energy]\cite[Theorem 7.2]{EKYY2}
 \label{bulkWigner}
Suppose that $H = (h_{ij})$ is a complex Hermitian (respectively, real symmetric) generalized Wigner matrix.
 Suppose that for some constants $\e>0$, $C>0$,  
\begin{equation} \label{4+e}
\E \left | \sqrt N  h_{ij}  \right | ^{4 + \e}  \;\leq\; C.
\end{equation}
Let $n \in \N$ and $O : \R^n \to \R$ be a test function (i.e. compactly supported and continuous).
Fix $|E_0|<2$ and  $\xi > 0$,
then with  $b_N=N^{-1 + \xi}$
we have
\begin{multline}\label{avg}
\lim_{N \to \infty} \int_{E_0 - b_N}^{E_0 + b_N} \frac{\rd E}{2 b_N} \int_{\bR^n} \rd \alpha_1 \cdots
 \rd \alpha_n\, O(\alpha_1, 
\dots, \alpha_n) 
\\ 
{}\times{}  \frac{1}{\varrho_{sc}(E)^n} \left ( {p_{N}^{(n)} - p_{{\rm G}, N}^{(n)}} \right ) 
 \left ( {E +
\frac{\alpha_1}{N\varrho_{sc}(E)}, \dots, E + \frac{\alpha_n}{N\varrho_{sc}(E)}}\right )  \;=\; 0\,.
\end{multline}
Here $\varrho_{sc}$ is the semicircle law
 defined in \eqref{def:sc}, $p_N^{(n)}$ is the $n$-point correlation function
 of the eigenvalue 
distribution of $H$ \eqref{def:corr}, and $p^{(n)}_{{\rm G},N}$ is the $n$-point correlation function
of an $N \times N$ GUE (respectively, GOE) matrix.
\end{theorem}
The additional rescaling in \eqref{avg} with $\varrho_{sc}(E)$ is not essential, it just
reflects the choice of variables under which the Gaussian correlation function is
given exactly by the sine kernel \eqref{sine} and not by some trivially rescaled version of it.

We remark that our method also provides an effective speed of convergence in \eqref{avg}.
We also point out that the condition \eqref{1.3C} can be  relaxed, see Corollary 8.3 \cite{EKY3}. For example, 
the lower bound can be changed to  $N^{-9/8+\e}$. Alternatively, 
under an additional symmetry condition on the law of the matrix elements, the upper
bound can be relaxed to $N^{-8/9-\e}$.

For the next result, we introduce the notation $\llbracket A, B\rrbracket : = \{ A, A+1, \ldots, B\}$
for any integers $A<B$. A relatively straightforward consequence of Theorem~\ref{bulkWigner} is the average gap universality:

\begin{corollary}[Gap universality with averaged label]
 \label{gapWigner} 
Let  $H$ be as in Theorem~\ref{bulkWigner} and $O$ be a test function of $n$ variables.
Fix small positive constants $\xi,\al>0$.
Then for 
 any integer $j_0 \in  \llbracket\al N, (1-\al)N\rrbracket$ 
we have
\be\label{EEOave}
  \lim_{N\to\infty}\frac{1}{2N^\xi}  \sum_{|j-j_0|\le N^\xi} 
\big[\E -\E^{G}\big] O\big( N(\la_j-\la_{j+1}), N(\la_j-\la_{j+2}), \ldots , N(\la_j - \la_{j+n})\big) =0.
\ee
Here  $\lambda_j$'s are the ordered eigenvalues. $\E$ and $\E^G$ denote the expectation with respect to the
Wigner ensemble $H$ and the Gaussian (GOE or GUE) ensemble, 
 respectively.
\end{corollary}
We remark that, similarly to the explicit formulas for
the correlation functions \eqref{sine}, for Gaussian (GOE or GUE) ensembles
 there are explicit expressions for  the gap distribution even without local averaging.
They are given in terms of  a Fredholm determinant of the corresponding kernel $K$,
see \cite{De1, DG1, Sch}.

Now we present our results for fixed energy:

\begin{theorem}[Universality at fixed energy]\cite{BEYY3}
 \label{bulkWignerfixed} Theorem~\ref{bulkWigner} holds 
under the same conditions
 without averaging, i.e.
for 
any $E$ with $|E|<2$ we have
\begin{multline}\label{avgfix}
\lim_{N \to \infty}  \int_{\bR^n} \rd \alpha_1 \cdots
 \rd \alpha_n\, O(\alpha_1, 
\dots, \alpha_n) 
\\ 
{}\times{}  \frac{1}{\varrho_{sc}(E)^n} \left ( {p_{N}^{(n)} - p_{{\rm G}, N}^{(n)}} \right ) 
 \left ( {E +
\frac{\alpha_1}{N\varrho_{sc}(E)}, \dots, E + \frac{\alpha_n}{N\varrho_{sc}(E)}}\right )  \;=\; 0\,.
\end{multline}
\end{theorem}

We remark that the fixed energy result \eqref{avgfix} for the $\beta=2$ (complex Hermitian) case was
already known before, see \cite{EPRSY, TV5} for special cases and \cite{EY} for the general case.
 The $\beta=2$ case is exceptional since
the Harish-Chandra/Itzykson/Zuber  identity allows one to compute
correlation functions for Wigner matrices with a tiny Gaussian component.
This method relies on an algebraic identity and cannot be generalized to
other symmetry classes.

\medskip

Finally, the gap universality with fixed label
asserts that \eqref{EEOave} holds without averaging.

\begin{theorem}  [Gap universality with fixed label]\cite[Theorem 2.2]{EYsinglegap}
 \label{thm:sg} 
Assuming subexponential decay of the matrix elements instead of \eqref{4+e},
Corollary~\ref{gapWigner} holds without averaging:
\be\label{EEO}
 \lim_{N\to\infty}  \big[\E -\E^G\big] O\Big( N(\la_j-\la_{j+1}), N(\la_j-\la_{j+2}), \ldots , N(\la_j - \la_{j+n})\Big) =0,
\ee
for any $j\in \llbracket \al N, (1-\al)N\rrbracket$ with a fixed $\al>0$.

More generally, for any $k, m \in \llbracket \al N, (1-\al)N\rrbracket$ we have
\begin{align} \label{EEO1}
   \lim_{N\to\infty} \Big| \E  O\Big( & (N\varrho_k) (\la_k-\la_{k+1}),   (N\varrho_k)(\la_k-\la_{k+2}), \ldots , (N\varrho_k)(\la_k - \la_{k+n})\Big) \\
& -\E^G O\Big( (N\varrho_m)(\la_m-\la_{m+1}), \ldots , (N\varrho_m)(\la_m - \la_{m+n})\Big)\Big|
   =0,  \non
\end{align}
where the local density $\varrho_k$ is defined by
$\varrho_k : =\varrho_{sc}(\gamma_k)$ with $\gamma_k$  from  \eqref{def:gamma}.
\end{theorem}
The second part  \eqref{EEO1} of this theorem asserts that the gap distribution is not only
independent of the specific Wigner ensemble, but it is also universal throughout
the bulk spectrum. This is the counterpart of the statement that the appropriately
rescaled correlation functions \eqref{avgfix} have a limit that is independent of $E$, see
\eqref{sine}.

Prior to our work, universality for a single gap  was only achieved in  the special 
case of the Gaussian unitary ensemble (GUE) in \cite{Taogap}, which statement  then  easily  implies
 the same results for  complex Hermitian Wigner  matrices satisfying  the  
  four moment matching condition.

\smallskip

\subsection{Log-gases}

In the case of  invariant ensembles, it is well-known that for 
 $V$ satisfying certain mild conditions the sequence
of one-point correlation functions, or
 densities, associated with  $\mu=\mu^{(N)}_{\beta,V}$ from \eqref{loggas} has a limit
as $N\to\infty$
and the limiting  equilibrium density  $\varrho_V(s)$
can be obtained as the unique minimizer  of the
functional
$$
I(\nu)=
\int_\bR V(t) \nu(t)\rd t-
\int_\bR \int_\bR \log|t-s| \nu(s) \nu(t) \rd t \rd s.
$$
We assume that 
$\varrho=\varrho_{V}$ is supported on a single compact interval, $[A,B]$ and $\varrho\in C^2(A,B)$. 
Moreover, we assume that  $V$ is {\it regular} in the sense that $\varrho$ is strictly positive on $(A, B)$
and vanishes as a square root at the endpoints, i.e.
\begin{align}\label{sqsing}
\varrho(t)&=s_A\sqrt{t-A}\left(1+O\left(t-A\right)\right),\ t\to A^+,
\end{align}
for some constant $s_A>0$ and a similar condition holds at the upper edge. 

It is known that these conditions are satisfied if, for example, $V$ is strictly convex.
In this case $\varrho_V$ satisfies the equation 
\be
   \frac{1}{2}V'(t) = \int_\bR \frac{\varrho_V(s)\rd s}{t-s}
\label{equilibrium}
\ee
for any $t\in(A,B)$.  For the Gaussian case, $V(x)=x^2/2$, the equilibrium density
is given by the semicircle law, $\varrho_V=\varrho_{sc}$, see \eqref{def:sc}.

The following result was proven in Corollary  2.2 of \cite{BEY}
for convex real analytic potential $V$, it was generalized in Theorem 1.2 of \cite{BEY2} 
for the non-convex case and further generalized for arbitrary $C^4$ potential in Theorem 2.5 of \cite{BEY3}.

\begin{theorem}[Universality with averaged energy] \label{bulkbeta}
 Assume
$V\in C^4(\bR)$, regular and
let $\beta> 0$. Consider the $\beta$-ensemble $\mu_V=\mu_{\beta, V}^{(N)}$ given in \eqref{loggas} 
with correlation functions
 $p_{V,N}^{(n)}$  defined analogously to \eqref{def:corr}.
 For the Gaussian case, $V(x) =x^2/2$, the correlation
functions are denoted by $p_{G,N}^{(n)}$.
Let $E_0\in (A,B)$ lie in the interior of the support of $\varrho$
 and
similarly let $E'_0\in (-2,2)$ be inside the support of $\varrho_{sc}$. 
Then for $b_N = N^{-1+\xi}$ with some $\xi>0$  we have
\begin{align} \label{bulkbetaeq}
\lim_{N \to \infty} \int  & \rd \alpha_1 \cdots \rd \alpha_n\, O(\alpha_1, 
\dots, \alpha_n) \\
\nonumber
\times \Bigg [ & 
   \int_{E_0 - b_N}^{E_0 + b_N} \frac{\rd E}{2 b_N}  \frac{1}{ \varrho (E)^n  }  p_{V,N}^{(n)}   \Big  ( E +
\frac{\alpha_1}{N\varrho(E)}, \dots,   E + \frac{\alpha_n}{N\varrho(E)}  \Big  ) \\ \nonumber
&
 -   \int_{E'_0 - b_N}^{E'_0 + b_N} \frac{\rd E'}{2b_N}  \frac{1}{\varrho_{sc}(E')^n} p_{{\rm G}, N}^{(n)}      \Big  ( E' +
\frac{\alpha_1}{N\varrho_{sc}(E')}, \dots,   E' + \frac{\alpha_n}{N\varrho_{sc}(E')}  \Big  ) \Bigg ]
 \;=\; 0\,,
\end{align}
 i.e. the correlation functions  of $\mu_{\beta, V}^{(N)}$ averaged around  $E_0$
 asymptotically coincide with those
of the Gaussian case. In particular,  they are independent of $E_0$.
\end{theorem}
Theorem~\ref{bulkbeta} immediately implies gap universality with averaged label, exactly in
the same way as Corollary~\ref{gapWigner} was deduced from Theorem~\ref{bulkWigner}; we refrain
from stating it explicitly.
The following two theorems show that these results hold without averaging.

\begin{theorem}[Universality at fixed energy]\cite{BEYY3} \label{bulkbetafixed}
Consider the setup of Theorem~\ref{bulkbeta}
and we additionally assume that $\beta\ge 1$. 
Then \eqref{bulkbetaeq} holds without averaging, i.e.
for any $E\in (A,B)$ and $E'\in (-2,2)$ we have
\begin{align} \label{bulkbetaeqfix}
\lim_{N \to \infty} \int  & \rd \alpha_1 \cdots \rd \alpha_n\, O(\alpha_1, 
\dots, \alpha_n) \\
\nonumber
\times \Bigg [ & 
   \frac{1}{ \varrho (E)^n  }  p_{V,N}^{(n)}   \Big  ( E +
\frac{\alpha_1}{N\varrho(E)}, \dots,   E + \frac{\alpha_n}{N\varrho(E)}  \Big  ) \\ \nonumber
&
 -     \frac{1}{\varrho_{sc}(E')^n} p_{{\rm G}, N}^{(n)}      \Big  ( E' +
\frac{\alpha_1}{N\varrho_{sc}(E')}, \dots,   E' + \frac{\alpha_n}{N\varrho_{sc}(E')}  \Big  ) \Bigg ]
 \;=\; 0\,.
\end{align}
\end{theorem}
Prior to our work and with a different method,
the same result was also proven in \cite{Sch2} for analytic potentials and for any $\beta>0$
even if the support of $\varrho$ has several intervals. An extension of the method  to $V\in C^5$ is anticipated
in \cite{Sch2}.

To formulate the result for the gap universality with a fixed label, 
we define the quantiles $\gamma_{j,V}$ of the density $\varrho_V$ by
\be\label{defgammagen}
   \frac{j}{N} = \int_A^{\gamma_{j,V}}\varrho_V(x) \rd x,
\ee
similarly  to \eqref{def:gamma}.
We set
\be\label{rhov}
   \varrho_j^V := \varrho_V(\gamma_{j,V}),\qquad
\mbox{and}\qquad \varrho_j : = \varrho_{sc} (\gamma_{j})
\ee
to be the limiting densities at the $j$-th quantiles.
Let $\E^{\mu_V}$ and $\E^G$ denote the expectation w.r.t. the measure $\mu_V$ and
its Gaussian counterpart for $V(\lambda)=\frac{1}{2}\lambda^2$.

\begin{theorem}  [Gap universality with fixed label]\cite[Theorem 2.3]{EYsinglegap}
 \label{thm:beta} Consider the setup of Theorem~\ref{bulkbeta}
and we also assume $\beta\ge1$.  Set some $\al>0$, then
\begin{align}\label{betaeq}
 \lim_{N\to\infty} \Bigg| \E^{\mu_V} O\Big( & (N\varrho_k^V) (\la_k-\la_{k+1}),  (N\varrho_k^V) (\la_{k}-\la_{k+2}), \ldots , 
(N\varrho_k^V) (\la_k - \la_{k+n})\Big) \\ \nonumber
  &  - \E^{\mu_{G}} 
  O\Big( (N\varrho_m) (\la_m-\la_{m+1}), 
 \ldots , 
(N\varrho_m)(\la_m - \la_{m+n})\Big)
 \Bigg| 
 =0
\end{align}
for any $k,m\in \llbracket \al N, (1-\al)N\rrbracket$. 
In particular, the distribution of the rescaled gaps
 w.r.t. $\mu_V$ 
does not depend on the index $k$ in the bulk.
\end{theorem}

We point out that Theorem~\ref{bulkbeta} holds for any $\beta>0$, but Theorems~\ref{bulkbetafixed} and
\ref{thm:beta}
require $\beta\ge 1$. This is only a technical restriction 
related to a certain condition in the De Giorgi-Nash-Moser regularity theory
that is the backbone of our proof.
Indeed, a year after our work was completed, an alternative proof of \eqref{betaeq} was given for any $\beta>0$
but with a higher regularity assumption on $V$ and with an additional hypothesis that
can be effectively checked only for convex $V$, see \cite{BFG}.

\subsection{Universalities at the edge}

We stated our results for the bulk of the spectrum. Similar results hold at the edge; in
this case the ``averaged'' results are meaningless.  For completeness, we give the
 universality results for both ensembles.

\begin{theorem}  [Universality at the edge for Wigner matrices]\cite{BEY3} \label{thm:wigneredge}
Let $H$ be a generalized
Wigner ensemble
with subexponentially decaying matrix elements.
Fix $n\in \N$, $\kappa<1/4$ and a test function $O$ of $n$ variables.
Then for any $\Lambda\subset\llbracket 1,N^\kappa\rrbracket$
with $|\Lambda|= n$, we have
$$
\left|
\big[ \E-\E^{G}\big]
O
\left(\left(N^{ 2/3\nc} j^{1/3}(\lambda_j-\gamma_j)\right)_{j\in \Lambda}\right)
\right|\leq  N^{-\chi},
$$
with some $\chi>0$,
where $\E^G$  is expectation w.r.t. the standard GOE or GUE  ensemble depending on the symmetry class of $H$
and $\gamma_j$'s are semicircle quantiles.
\end{theorem}

Edge universality for Wigner matrices was first proved in \cite{Sos1999}
assuming symmetry of the distribution of the matrix elements and finiteness
of all their moments. The symmetry condition was completely eliminated \cite{EYYrigi}
and the optimal moment condition  was obtained in \cite{LY}.
All these works heavily rely on the fact that the variances of
the matrix elements are identical. The main point of Theorem~\ref{thm:wigneredge}
is to consider generalized Wigner
matrices, i.e., matrices with non-constant variances.
In fact, it was shown in \cite{EYYrigi}  that the edge statistics for any generalized Wigner
matrix are universal in the sense that they coincide with those of a generalized Gaussian Wigner
matrix with the same variances, but it was not shown that the statistics
are independent of the variances themselves.  Theorem~\ref{thm:wigneredge}
provides this missing step and thus it proves the edge universality
in the broadest sense.

\begin{theorem}[Universality at the edge for log-gases]\label{thm:betaedge}\cite{BEY3}
Let  $\beta\ge 1$ and $V$ (resp. $\widetilde V$) be in $C^4(\bR)$,
regular such that the equilibrium density $\varrho_V$ (resp. $\varrho_{\widetilde V}$) is
supported on a single interval $[A,B]$ (resp. $[\widetilde A,\widetilde B]$).
 Without loss of generality
  we  assume  that for both densities $(\ref{sqsing})$
 holds with $A=0$ and with the same constant $s_A$.
Fix $n\in \N$,  $ \kappa <  2/5$.
Then for any  $\Lambda\subset\llbracket 1,N^\kappa\rrbracket$
with $|\Lambda|= n$, we have
\be\label{mm}
\left|
(\E^{\mu_V}-\E^{\mu_{\widetilde V}})
O
\left(\left(N^{ 2/3} j^{1/3}(\lambda_j-\gamma_j)\right)_{j\in \Lambda}\right)
\right|\leq  N^{-\chi}
\ee
with some $\chi>0$.
Here $\gamma_j$ are the quantiles w.r.t. the density $\varrho_V$ \eqref{defgammagen}.
\end{theorem}

The first results on edge universality for invariant ensembles concerned the classical values of $\beta=1,2,4$. 
The case 
$\beta=2$ and real analytic $V$ was  solved in \cite{DKMVZ1, DGedge}.
The $\beta=1,4$ cases are  considerably harder than $\beta=2$. For  $\beta=1,4$
universality was first proved  for polynomial potentials in  \cite{DGedge},
then for the real analytic case for $\beta=1$  in \cite{PSedge, Schedge}, which also give an alternative proof for $\beta=2$. 
Finally, independently of our work with a completely different method, edge universality for any
$\beta>0$ 
and convex  polynomial $V$ was recently proved in \cite{KriRidVir2013}.

\section{Outline of the proof strategy}

\subsection{``Averaged'' results: Dyson Brownian motion}\label{sec:DBM}

The proof  of Theorem \ref{bulkWigner} 
follows a  three-step  strategy that was first introduced in \cite{EPRSY}
and further developed in \cite{ESY4}.

\medskip 
\noindent
{\it Step 1.}  {\it Local semicircle law and rigidity of eigenvalues.}
The main tool is the resolvent of $H$ at a spectral parameter $z=E+i\eta$
with $\eta\gg 1/N$;
$$
   m_N(z): = \frac{1}{N}\tr \frac{1}{H-z} = \frac{1}{N}\sum_{j} \frac{1}{\lambda_j-z},
$$
which is of the form of \eqref{sclawlocal} with $O(x)= (x-i)^{-1}$. Using the Schur
decomposition formula we may write
$$
  m_N(z) = \frac{1}{N} \sum_{j=1}^N \frac{1}{h_{jj} - z - \sum_{a,b\ne j} h_{ja} 
 G_{ab}^{(j)}(z)
 h_{bj}},
$$
where $G^{(j)}(z) = (H^{(j)}-z)^{-1}$ is the resolvent of the $(N-1)\times (N-1)$ minor $H^{(j)}$ of $H$ after removing the
$j$-th row and column. Since $G_{ab}^{(j)}(z)$ and $h_{ja} h_{bj}$ are independent, we may use 
concentration results to replace the double sum in the denominator by its 
expectation over the matrix elements in the $j$-th row and column. Neglecting the fluctuation, we recover  $m_N^{(j)}(z)$,
the normalized trace of the resolvent of $H^{(j)}$. Since $m_N^{(j)}(z)$ and $ m_N(z)$ are close, we obtain
the following {\it self-consistent equation} 
\be\label{sce}
  m_N(z) = -\frac{1}{z+ m_N(z)} + \mbox{error}.
\ee
If the error is neglected, then the solution of the resulting quadratic equation is
exactly the Stieltjes transform
$$
  m_{sc}(z): = \int_\bR \frac{1}{x-z} \varrho_{sc}(x) \rd x
$$
of the Wigner semicircle law $\varrho_{sc}(x)$. This allows us to conclude that $m_N(z)$ is close to $m_{sc}$,
and a careful analysis yields
\be\label{lsc}
   | m_N(z) - m_{sc}(z)|\lesssim \frac{1}{N\eta}.
\ee
This is the local semicircle law in resolvent form, from which
the limit of \eqref{sclawlocal}
and the rigidity property \eqref{rig} can be concluded.

\medskip 
\noindent
{\it Step 2.}
{\it Universality for
Gaussian divisible ensembles:}  The Gaussian divisible ensembles are matrices of the form 
$$
 H_t= e^{-t/2} H+ \sqrt {1 - e^{-t}} U,
$$
 where $H$ is a Wigner matrix and
$U$ is an independent GUE/GOE matrix. The parametrization of  $H_t$ reflects that, in the sense of
distribution, it is most conveniently obtained
by an Ornstein-Uhlenbeck  process:
\be\label{OU}
   \rd H_t = \frac{1}{\sqrt{N}} \rd B_t - \frac{1}{2} H_t \rd t,
\ee
where $B_t$ is a matrix-valued Brownian motion of the appropriate symmetry class. 
Dyson observed \cite{DyB} that the corresponding process $\bla_t$ of the eigenvalues
of $H_t$ remarkably satisfies a system of stochastic
differential equations (SDE), called the {\it Dyson Brownian Motion (DBM)}:
\be
    \rd \lambda_j = \frac{1}{\sqrt{N}}\rd B_j +\Big[ - \frac{1}{2}\lambda_j+ \frac{1}{N}
\sum_{k\ne j} \frac{1}{\lambda_k-\lambda_j} 
 \Big]\rd t,
\label{dbm}
\ee
written for $\beta=2$, where $B_j$'s are independent standard  real Brownian motions.
The key idea is to study the relaxation of the flow \eqref{dbm} to its equilibrium
measure which is the distribution of the GUE eigenvalues. It turns out that,
 tested against 
observables involving only {\it differences of eigenvalues}, the convergence is extremely fast.
Combined with the rigidity bound 
that guarantees a strong apriori control on the initial state, we obtain that the 
gap statistics are already in local equilibrium (hence universal) after a very short time $t=N^{-1+\e}$,
see \cite{ESY4, ESYY}.

This method substantially improves Johansson's result \cite{J} which showed universality
only with a substantial Gaussian component (essentially for $t>0$ independent of $N$)
and only for the $\beta=2$ symmetry class. In fact, the first restriction can be
relaxed by using our optimal rigidity bound \cite{EPRSY, EY}, but the second one cannot be removed since the
proof relies on the  Harish-Chandra/Itzykson/Zuber formula. The analysis of the DBM is
much more robust, in particular it applies to any symmetry class. However, it yields
only an averaged result \eqref{EEOave} (from which \eqref{avg} can be deduced), while \cite{EPRSY, EY}
gives the fixed energy results \eqref{avgfix} but only for $\beta=2$.

\medskip 
\noindent
{\it Step 3.}  {\it Approximation by Gaussian divisible ensembles:}
 It is a simple 
 density argument  in the space
of matrix ensembles which 
shows that for any probability distribution of the matrix
 elements there exists a Gaussian divisible distribution
with a small Gaussian component, as in Step 2, such 
that the two  associated Wigner ensembles
have asymptotically identical local eigenvalue statistics. 
 The first implementation  of this approximation scheme was via  a reverse heat flow 
argument  \cite{EPRSY}; it was later replaced by  the 
  {\it Green function comparison theorem} \cite{EYY}
motivated by the four moment matching condition of \cite{TV}. 
This comparison argument is very robust: it works even without averaging and
for arbitrary observables not only for those of difference type.

\medskip

\noindent
The proof of Theorem \ref{bulkbeta} follows a somewhat similar path 
but with essential differences. Rigidity estimates still hold on the
smallest scale, but their derivation cannot use resolvents
since there is no matrix behind a general log-gas.
Instead of \eqref{sce} we use the loop equation
from \cite{Joh} or \cite{Sch}, but 
extended to smooth potentials.  There is no analogue
of the Gaussian divisible ensemble for log-gases, but an enhanced version 
of the DBM underlying the invariant measure $\mu_V$ can 
still be analyzed. 

In summary, the DBM plays the fundamental role
behind the ``averaged'' universality result for both models.

\subsection{``Fixed'' results: H\"older regularity and homogenization} For definiteness, we will present some ideas
to prove Theorem~\ref{thm:beta}, the proof of Theorems~\ref{bulkWignerfixed}, \ref{thm:sg}, \ref{bulkbetafixed}
and the results at the edge require additional steps.

\medskip

{\it Step 1. Comparison of local Gibbs measures.} The basic mechanism for
universality is that the microscopic structure of the measure $\mu_V$ defined in \eqref{pinv}
is  insensitive of the potential $V$, 
it is essentially determined by the 
Vandermonde determinant, i.e. the log-interaction in \eqref{loggas}.
In the first step we localize the problem by freezing (conditioning on)
all particles at a distance $1\ll K \ll N$ away from the fixed 
index $j$ of the gap
$\lambda_j-\lambda_{j+1}$ we want to study. Thus the corresponding local Gibbs measure is defined 
on an interval $I= [j-K, j+K]$ and it still
retains the Vandermonde structure. On this mesoscopic scale the potential is locally constant, hence 
its effect is trivial, so the key question is to show that 
$\lambda_j-\lambda_{j+1}$ is largely insensitive to  the boundary effects
we just introduced by localization. This is a question about the 
long range correlation structure
of the Gibbs measure. 

The main difficulty is that the log-gas is a strongly correlated
system in  contrast to the
customary setup in statistical physics where correlations often
decay very fast. In fact, the covariance between two points decays only logarithmically
\be\label{corrr}
   \frac{ \langle \la_i ; \la_j\rangle}{\sqrt{  \langle \la_i ; \la_i\rangle \langle \la_j ; \la_j\rangle}}\sim \frac{1}{\log |i-j|}, 
 \qquad 1\ll |i-j|\ll N.
\ee
One key observation is that the correlation decay
 between a {\it gap} $\lambda_i-\lambda_{i+1}$ and a 
{\it point} $\lambda_j$ is faster, it is
$|i-j|^{-1}$,  practically the discrete derivative of \eqref{corrr}.

\medskip

{\it Step 2. Random walk representation of the covariance.}
In a more general setup, consider  a
Gibbs measure $\om(\rd\bx) = e^{-\beta \cH(\bx)}\rd \bx$ on finitely many points labelled by $I$  and with a strictly
convex Hamiltonian, $\cH''(\bx)\ge c>0$.  Then the covariance w.r.t. $\om$ can be expressed as
\be\label{rwr}
  \langle F(\bx); G(\bx) \rangle_\om = \frac{1}{2}\int_0^\infty
   \rd s \int \rd\om(\bx)\E_{\bx} \big[\;  \nabla G(\bx(s))
\cdot\cU(s, \bx(\cdot)) \nabla F(\bx) \big],
\ee
see \cite{HS, NS}.
Here $\E_\bx$ is the expectation for the (random) paths $\bx(\cdot)$ starting
from $\bx(0)=\bx$ and solving the canonical SDE for the measure $\om$:
\be\label{dbmgen}
    \rd \bx(s) =\rd {\bf B}(s) - \beta \nabla\cH ( \bx(s))\rd s,
\ee
and $\cU(s) =\cU(s, \bx(\cdot))$ is the fundamental solution to the linear
system of equations
\be\label{matrixeq}
   \partial_s \cU(s) =- \cU(s) \cA(s), \qquad \cA(s): = \beta \cH''(\bx(s))
\ee
with $\cU(0)=\mbox{Id}$. 
Notice that the coefficient matrix $\cA(s)$, and thus the fundamental
solution, depend on the random path. The SDE \eqref{dbmgen} is the generalization
of the DBM, \eqref{dbm}. Formula \eqref{rwr} turns the problem of
computing the covariance $ \langle F; G \rangle$ into a time-dependent question
to understand the fundamental solution $\cU$ of the parabolic equation \eqref{matrixeq}.

In particular, if $G$ is a function of a single  gap, $G(\bx)= O(x_j-x_{j+1})$
with some fixed $j$, and $F$ represents the
boundary effects, then \eqref{rwr} becomes
 \be\label{rwr1}
\frac{1}{2}\int_0^\infty
   \rd s \int \rd\om(\bx)\sum_{i\in I} 
  \E_\bx \Big[\; O'(x_j-x_{j+1}) \big(
\cU_{i,j}(s) -\cU_{i,j+1}(s)\big)  \pt_i F(\bx) \Big].
\ee
The key technical step is to show that for a typical
path $\bx(\cdot)$ the solution $\cU(s)$ is H\"older-regular in a sense that
$\cU_{i, j}(s)-\cU_{i, j+1}(s)$ is small if $j$ is away from the boundary of $I$ and
$s$ is not too small.

\medskip
{\it Step 3. H\"older-regularity of the solution to \eqref{matrixeq}.}
For any fixed realization of the path $\bx(\cdot)$, we will view 
the equation \eqref{matrixeq} as a finite dimensional version of
a parabolic equation. The coefficient matrix,  the Hessian of
the local Gibbs measure, is computed explicitly. It can be
 written as
$\cA = \cB +\cW$, where $\cW\ge 0$ is diagonal, $\cB$ is a symmetric matrix
with quadratic form
$$
    \langle \bu, \cB(s) \bu\rangle = \frac{1}{2} \sum_{i,j\in I} B_{ij}(s) (u_i-u_j)^2,
  \qquad B_{ij}(s):= \frac{\beta}{(x_i(s)-x_j(s))^2}.
$$
After rescaling the problem so that the gap is of order one,
for a typical path and large $i-j$ we have 
\be\label{cbs}
   B_{ij}(s) \sim \frac{1}{(i-j)^2}
\ee
by rigidity. We also have a lower bound for any $i\ne j$
\be\label{lowcbs}
  B_{ij}(s) \gtrsim \frac{1}{(i-j)^2},
\ee
at least with a very high probability.
If a matching upper bound were true for any $i\ne j$, then \eqref{matrixeq}
would be the discrete  analogue of the general  equation
\be
\partial_t u(t, x) = \int K(t, x, y) [u(t, y) - u(t, x) ] \rd y, \qquad t>0, \quad x, y\in \R^d 
\ee
  considered by  Caffarelli-Chan-Vasseur  in \cite{C}, where the 
 kernel $K$ is symmetric and has a specific short distance singularity
\be\label{Kxyt1}
C_1  |x-y|^{-d-s} \le K(t, x, y) \le C_2  |x-y|^{-d-s}
\ee
for some $s\in (0,2)$ and positive constants $C_1, C_2$. 
Roughly speaking, the integral operator  $K$ corresponds to the behavior 
of the operator $|p|^{s}$,  where $p=-i\nabla$.  The main result of \cite{C} asserts
that for any $t_0>0$, the solution $u(t,x)$ is $\e$-H\"older continuous,
 $u\in C^\e( (t_0, \infty), \R^d)$,
for some positive exponent $\e$  that depends only on $t_0$, $C_1$, $C_2$.
This is a version of the celebrated De Giorgi-Nash-Moser regularity result
for a non-local operator.

Our equation \eqref{matrixeq} is of this type with $d=s=1$, but it
is discrete and in a finite interval $I$ with a potential term. 
The key difference, however, is
that the coefficient  $B_{ij}(t)$  
can be singular in the sense that $B_{ij}(t) |i-j|^2$ is not uniformly bounded
 when $i,j$ are close to each other.
 Thus the analogue of the uniform upper bound  \eqref{Kxyt1} does not even hold for a fixed $t$. 
We first need to regularize the singularity of $B_{ij}$ on a very tiny scale. 
Even after  that we  can control  the regularized $B^{\tiny \mbox{reg}}_{ij}$   \normalcolor only in a certain
 average sense:
\be\label{Kass1}
   \sup_{0 \le s \le \si}\sup_{0 \le M\le K} \frac{1}{ 1+ s} \int_0^s \frac{1}{M}
 \sum_{i\in I\, : \, |i-Z| \le M}
 B_{i,i+1}^{\tiny \mbox{reg} \normalcolor}(s)
 \rd s \le CK^{\rho}
\ee
with high probability, for some small exponent $\rho$ 
and for any fixed $Z$ away from the edges of $I$.
 This estimate essentially says that 
the space-time maximal function of $B_{i, i+1}^{ \tiny \mbox{reg} \normalcolor}(t)$   at a fixed space-time point $(Z,0)$ 
is bounded by $K^\rho$. 
Our main generalization of the result in \cite{C} is to show
that the weak upper bound \eqref{Kass1} at a few space-time points together with \eqref{cbs} and \eqref{lowcbs}
(holding up to a factor $K^\xi$)
are sufficient for proving a discrete version of the H\"older continuity at the point $(Z,0)$.
More precisely,  there exists an $\e>0$ such that 
for any fixed $1\ll \si\ll K$
 the solution to \eqref{matrixeq} satisfies
\be\label{holds}
   \sup_{|j-Z|+ |j'-Z|\le \si^{1-\al}}  | \cU_{i, j}(\si) - \cU_{i, j'}(\si)| \le C K^\xi
\si^{-1- \e\al}
\ee
with any $\al\in[0,1/3]$  if we can 
guarantee that $\rho$ and $\xi$ are sufficiently small. 
The exponent $\e$  plays the
role of the H\"older regularity exponent.
Notice that 
 $\cU_{i,j}(\si)$ decays as $\si^{-1}$, hence \eqref{holds} provides an
additional decay for the discrete derivative. In particular, this
guarantees that the $\rd s$ integration in \eqref{rwr1} is finite.
With several further technical steps, this proves Theorem~\ref{thm:beta}.

\medskip

{\it Step 4. Homogenization.}   The proofs of Theorems~\ref{thm:sg}, \ref{bulkbetafixed}
require an additional information about the fundamental solution of \eqref{matrixeq}.
Since in the $|i-j|\gg 1$ regime we have $B_{ij}(s) \sim |i-j|^{-2}$, it is reasonable
to expect that the large time and large scale behavior of $\cU$ is given by the 
\be
  \cU_{ij}(t) \approx \Big(e^{-t|p|}\Big)_{ij} = \frac{t}{t^2 + (i-j)^2}, \quad |i-j|\gg 1,\;\;  t\gg 1,
\label{homog}
\ee
where we computed the heat kernel of $|p|=\sqrt{-\Delta}$ explicitly. This result, combined with a coupling
argument, yields that
\be\label{keyf}
     \lambda_i(t) - \wt \lambda_i (t) =  \Big(e^{-t|p|} \bla(0) \Big)_i  - \Big(e^{-t|p|} \wt\bla(0)\Big)_i  + \mbox{error},
\ee
where $\bla$ and $\wt\bla$ are two solutions of the SDE \eqref{dbmgen} with the same Brownian motion ${\bf B}(s)$
but with two different initial conditions. In the applications, $\bla(0)$ will be GUE/GOE eigenvalues 
and $\wt\bla(0)$ will be the eigenvalues of a general Wigner matrix.
Formula \eqref{keyf} allows us to express a single Wigner eigenvalue $ \lambda_i(t)$
in terms of the corresponding Gaussian eigenvalue $ \wt \lambda_i (t)$ and in terms of averaged quantities involving many
eigenvalues. Since averaged quantities can be computed much easier and Gaussian computations can be
performed by explicit formulas, we obtain nontrivial information about $ \lambda_i(t)$. Finally, 
approximation ideas similar to Step 3. 
in Section~\ref{sec:DBM} can relate general Wigner eigenvalues to Wigner eigenvalues
 with some Gaussian component such as 
 $ \lambda_i(t)$. In particular, these ideas can prove 
the logarithmic correlation decay \eqref{corrr} for any Wigner matrix.

In summary, the detailed analysis of the parabolic equation \eqref{matrixeq} with
singular coefficients given by the Dyson Brownian motion play the crucial role behind all
``fixed'' universality results for both Wigner matrices and  log-gases.

\medskip

{\it Acknowledgement.}
 Most results in this paper were  obtained in collaboration
with Horng-Tzer Yau, Benjamin Schlein, Jun Yin, Antti Knowles and 
Paul Bourgade and in some work, also with
Jose Ramirez and Sandrine P\'eche.
 This article  reports the joint progress 
with these authors.


\end{document}